\documentclass[a4paper,12pt]{article}
\usepackage{amsmath, amssymb, amsthm}

\theoremstyle{plain}

\newtheorem{thm}{Theorem}[section]
\newtheorem{prop}[thm]{Proposition}
\newtheorem{lemma}[thm]{Lemma}

\newtheorem{corollary}[thm]{Corollary}

\theoremstyle{definition}

\begin{document}

\def\proof{\paragraph{Proof.}}
\def\endproof{\hfill$\square$}
\def\noproof{\hfill$\square$}

\def\Z{{\mathbb Z}}\def\N{{\mathbb N}} \def\C{{\mathbb C}}
\def\Q{{\mathbb Q}}\def\R{{\mathbb R}} \def\E{{\mathbb E}}
\def\P{{\mathbb P}}

\def\lg{{\rm lg}}\def\Id{{\rm Id}}\def\GG{{\cal G}}
\def\AA{{\cal A}}\def\cd{{\rm cd}}\def\mf{{\rm mf}}
\def\rkAb{{\rm rkAb}}\def\rkZ{{\rm rkZ}}\def\Ab{{\rm Ab}}
\def\HH{{\cal H}}\def\Ker{{\rm Ker}}

%%%%%%%%%%%%%%%%%%%%%%%%%%%%%%%%%%%%%%%%%%%%%%%%%%%%%%%%%%%%%
\title{\bf{Artin groups of spherical type up to isomorphism}}
 
\author{
\textsc{Luis Paris}}

\date{\today}

\maketitle

\begin{abstract} 
We prove that two Artin groups of spherical type are isomorphic if and only if their defining 
Coxeter graphs are the same.
\end{abstract}

\noindent
{\bf AMS Subject Classification:} Primary 20F36. 

%%%%%%%%%%%%%%%%%%%%%%%%%%%%%%%%%%%%%%%%%%%%%%%%%%%%%%%%%%%%%%%%%%%%%%%%%%%%%%%%%%%%%%%%%
\section{Introduction}

Let $S$ be a finite set. Recall that a {\it Coxeter matrix} over $S$ is a matrix 
$M=(m_{s\,t})_{s,t \in S}$ indexed by the elements of $S$ such that $m_{s\,s}=1$ for all $s \in 
S$, and $m_{s\,t}=m_{t\,s} \in \{2, 3, 4, \dots, +\infty\}$ for all $s,t \in S$, $s \neq t$. A 
Coxeter matrix $M=(m_{s\,t})$ is usually represented by its {\it Coxeter graph}, $\Gamma$, which 
is defined as follows. The set of vertices of $\Gamma$ is $S$, two vertices $s,t$ are joined by 
an edge if $m_{s\,t}\ge 3$, and this edge is labelled by $m_{s\,t}$ if $m_{s\,t} \ge 4$. For $s,t 
\in S$ and $m \in \Z_{\ge 2}$, we denote by $w(s,t:m)$ the word $sts \dots$ of length $m$. The {\it Artin 
group} associated to $\Gamma$ is defined to be the group $G=G_\Gamma$ presented by
$$
G= \langle S\ |\ w(s,t:m_{s\,t}) = w(t,s:m_{s\,t})\text{ for } s,t \in S,\ s \neq t\text{ and } 
m_{s\,t} < +\infty \rangle\,.
$$
The {\it Coxeter group} $W=W_\Gamma$ associated to $\Gamma$ is the quotient of $G$ by the 
relations $s^2=1$, $s\in S$. We say that $\Gamma$ (or $G$) is of {\it spherical type} if $W$ is 
finite, that $\Gamma$ (or $G$) is {\it right-angled} if $m_{s\,t} \in \{2, +\infty \}$ for all 
$s,t \in S$, $s \neq t$, and that $G$ (or $W$) is {\it irreducible} if $\Gamma$ is connected. The 
number $n=|S|$ is called the {\it rank} of $G$ (or of $W$). 

One of the main question in the subject is the classification of Artin groups up to isomorphism 
(see \cite{Bes2}, Question 2.14). This problem is far from being completely solved as Artin groups 
are poorly understood in general. 
For example, we do not know whether all Artin groups are torsion free, and we do not know 
any general solution to the word problem for these groups. The only known results concerning this 
classification question are contained in a work by Brady, McCammond, M\"uhlherr, and Neumann \cite{BMMN}, 
where the authors determine a sort of transformation on Coxeter graphs which does not change the 
isomorphism class of the associated Artin groups, and a work by Droms \cite{Dro}, where it is proved 
that, if $\Gamma$ and $\Omega$ are two right-angled Coxeter graphs whose associated Artin groups 
are isomorphic, then $\Gamma= \Omega$. Notice that an Artin group is biorderable if and only if 
it is right-angled, hence a consequence of Droms' result is that, if $\Gamma$ is a right-angled 
Coxeter graph and $\Omega$ is any Coxeter graph, and if the Artin groups associated to $\Gamma$ 
and $\Omega$ are isomorphic, then $\Gamma=\Omega$. The fact that right-angled Artin groups are 
biorderable is proved in \cite{DuTh}. In order to show that the remainig Artin groups are not 
biorderable, one has only to observe that, if $2<m_{s\,t}< +\infty$, then $(st)^{m_{s\,t}} = 
(ts)^{m_{s\,t}}$ and $st \neq ts$, and that, in a biorderable group, two distinct elements cannot 
have a common $m$-th power for a fixed $m$.

In this paper we answer the classification question in the restricted framework of spherical type
Artin groups. More precisely, we prove the following.

\begin{thm}
Let $\Gamma$ and $\Omega$ be two spherical type Coxeter graphs, and let $G$ and $H$ be the Artin 
groups associated to $\Gamma$ and $\Omega$, respectively. If $G$ is isomorphic to $H$, then 
$\Gamma=\Omega$.
\end{thm}

\noindent
{\bf Remark.} I do not know whether a non spherical type Artin group can be isomorphic to a 
spherical type Artin group.

\bigskip
Artin groups were first introduced by Tits \cite{Tit2} as extensions of Coxeter groups. Later, 
Brieskorn \cite{Bri1} gave a topological interpretation of the Artin groups of spherical type in 
terms of complements of discriminantal varieties.
Define a (real) {\it reflection group} of rank $n$ to be a finite subgroup $W$ of $GL(n,\R)$ 
generated by reflections. Such a group is called {\it essential} if there is no non-trivial 
subspace of $\R^n$ on which $W$ acts trivially. Let $\AA$ be the set of reflecting hyperplanes of 
$W$, and, for $H \in \AA$, let $H_\C$ denote the complexification of $H$, {\it i.e.} the complex 
hyperplane in $\C^n$ defined by the same equation as $H$. Then $W$ acts freely on $M(W)= \C^n 
\setminus \cup_{H \in \AA} H_\C$, and, by \cite{Che}, $N(W)= M(W)/W$ is the complement in $\C^n$ 
of an algebraic variety, $D(W)$, called {\it discriminantal variety} of type $W$. Now, take 
a spherical type Coxeter graph $\Gamma$, and consider the associated Coxeter group 
$W=W_\Gamma$. By \cite{Cox}, the group $W$ can be represented as an essential reflection group in 
$GL(n,\R)$, where $n=|S|$ is the rank of $W$, and, conversely, any essential reflection group of 
rank $n$ can be uniquely obtained in this way. By \cite{Bri1}, $\pi_1(N(W))$ is the Artin group 
$G=G_\Gamma$ associated to $\Gamma$.

So, a consequence of Theorem 1.1 is that $\pi_1(N(W))= \pi_1( \C^n \setminus D(W))$ completely 
determines the reflection group $W$ as well as the discriminantal variety $D(W)$.

Since the work of Brieskorn and Saito \cite{BrSa} and that of Deligne \cite{Del}, the 
combinatorial theory of spherical type Artin groups has been well studied. In particular, these 
groups are know to be biautomatic (see \cite{Cha2}, \cite{Cha}), and torsion-free. This last 
result is a direct consequence of \cite{Del} and \cite{Bri1}, it is explicitely proved in 
\cite{Deh}, and it shall be of importance in the remainder of the paper.

The first step in the proof of Theorem 1.1 consists of calculating some invariants for 
spherical type Artin groups (see Section 3). It actually happens that these invariants separate 
the irreducible Artin groups of spherical type (see Proposition 5.1). Afterwards, 
for a given isomorphism $\varphi: G \to H$ between spherical type Artin groups, we 
show that, up to some details, $\varphi$ sends each irreducible component of $G$ injectively into 
a unique irreducible component of $H$, and that both components have the same invariants. In 
order to do that, we first need to show that an irreducible Artin group $G$ cannot be decomposed 
as a product of two subgroups which commute, unless one of these subgroups lies in the center of 
$G$ (see Section 4).

From now on, $\Gamma$ denotes a spherical type Coxeter graph, $G$ denotes its associated Artin 
group, and $W$ denotes its associated Coxeter group.

\bigskip\noindent
{\bf Acknowledgments.} The idea of looking at centralizers of ``good'' elements in the 
proof of Proposition 4.2 is a suggestion of Benson Farb. I am grateful to him for this clever idea 
as well as for all his useful conversations. I am also grateful to Jean Michel who pointed out to 
me his work with Michel Brou\'e, and to John Crisp for so many discussions on everything 
concerning this paper.

%%%%%%%%%%%%%%%%%%%%%%%%%%%%%%%%%%%%%%%%%%%%%%%%%%%%%%%%%%%%%%%%%%%%%%%%%%%%%%%%%%%%%%%%%
\section{Preliminaries}

We recall in this section some well-known results on Coxeter groups and Artin groups.

For a subset $X$ of $S$, we denote by $W_X$ the subgroup of $W$ generated by $X$, and by $G_X$ the 
subgroup of $G$ generated by $X$. Let $\Gamma_X$ be the full Coxeter subgraph of $\Gamma$ whose 
vertex set is $X$. Then $W_X$ is the Coxeter group associated to $\Gamma_X$ (see \cite{Bou}), and $G_X$ 
is the Artin group associated to $\Gamma_X$ (see \cite{Lek} and \cite{Par1}). The subgroup $W_X$ is 
called {\it standard parabolic subgroup} of $W$, and $G_X$ is called {\it standard parabolic 
subgroup} of $G$. 

For $w \in W$, we denote by $\lg (w)$ the word length of $w$ with respect to $S$. The group $W$ has 
a unique element of maximal length, $w_0$, which satisfies $w_0^2=1$ and $w_0 S w_0 =S$, and whose 
length is $m_1 + \dots + m_n$, where $m_1, m_2, \dots, m_n$ are the exponents of $W$. 

The connected spherical Coxeter graphs are exactly the graphs $A_n$ ($n \ge 1$), $B_n$ ($n\ge 2$), 
$D_n$ ($n\ge 4$), $E_6$, $E_7$, $E_8$, $F_4$, $H_3$, $H_4$, $I_2(p)$ ($p \ge 5$) represented in 
\cite{Bou}, Ch. IV, \S 4, Thm. 1. (Here we use the notation $I_2(6)$ for the Coxeter graph $G_2$. 
We may also use the notation $I_2(3)$ for $A_2$, and $I_2(4)$ for $B_2$.)

Let $F:G \to W$ be the natural epimorphism which sends $s$ to $s$ for all $s \in S$. This 
epimorphism has a natural set-section $T: W \to G$ defined as follows. Let $w \in W$, and let 
$w=s_1s_2 \dots s_l$ be a reduced expression of $w$ ({\it i.e.} $l=\lg(w)$). Then $T(w)=s_1s_2 
\dots s_l \in G$. By Tits' solution to the word problem for Coxeter groups \cite{Tit}, the 
definition of $T(w)$ does not depend on the choice of the reduced expression.

Define the {\it Artin monoid} associated to $\Gamma$ to be the (abstract) monoid $G^+$ presented by
$$
G^+= \langle S\ |\ w(s,t: m_{s\,t}) = w(t,s: m_{s\,t}) \text{ for } s\neq t \text{ and } m_{s\,t} 
<+\infty \rangle^+\,.
$$
By \cite{BrSa}, the natural homomorphism $G^+ \to G$ which sends $s$ to $s$ for all $s \in S$ is 
injective. Note that this fact is always true, even if $\Gamma$ is not assumed to be of spherical 
type (see \cite{Par2}).

The {\it fundamental element} of $G$ is defined to be $\Delta= T(w_0)$, where $w_0$ denotes the 
element of $W$ of maximal length. For $X \subset S$,
We denote by $w_X$ the element of $W_X$ of maximal length,
and by $\Delta_X=T(w_X)$ the fundamental element of $G_X$.

The defining relations of $G^+$ being homogeneous, we can define two partial orders $\le_L$ and 
$\le_R$ on $G^+$ as follows.

\smallskip
$\bullet$ We set $a \le_L b$ if there exists $c \in G^+$ such that $b=ac$.

\smallskip
$\bullet$ We set $a \le_R b$ if there exists $c \in G^+$ such that $b=ca$.

\smallskip
Now, the following two propositions are a mixture of several well-known results from \cite{BrSa} 
and \cite{Del}.

\begin{prop}
(1) $G^+$ is cancellative.

\smallskip
(2) $(G^+, \le_L)$ and $(G^+, \le_R)$ are lattices.

\smallskip
(3) $\{a \in G^+; a\le_L \Delta\} = \{a \in G^+; a\le_R \Delta \} = T(W)$.
\noproof
\end{prop}

Note that the fact that $G^+$ is cancellative is true even if $\Gamma$ is not of spherical 
type (see \cite{Mic}). The elements of $T(W)$ are called {\it simple elements}. We shall denote 
the lattice operations of $(G^+, \le_L)$ by $\vee_L$ and $\wedge_L$,  
and the lattice operations of $(G^+, \le_R)$ by $\vee_R$ and $\wedge_R$.

Define the {\it quasi-center} of $G$ to be the subgroup 
$QZ(G)=\{ a \in G; aSa^{-1} = S \}$.

\begin{prop}
Assume $\Gamma$ to be connected.

\smallskip
(1) For $X \subset S$ we have
$$
\vee_L \{s; s\in X\} = \vee_R \{s; s\in X\} = \Delta_X\,.
$$
In particular,
$$
\vee_L \{s; s\in S\} = \vee_R \{s; s\in S\} = \Delta\,.
$$

\smallskip
(2) There exists a permutation $\mu: S \to S$ such that $\mu^2= \Id$ and $\Delta s = \mu(s) 
\Delta$ for all $s \in S$.

\smallskip
(3) The quasi-center $QZ(G)$ of $G$ is an infinite cyclic subgroup generated by $\Delta$.

\smallskip
(4) The center $Z(G)$ of $G$ is an infinite cyclic subgroup of $G$ generated either by $\delta= 
\Delta$ if $\mu=\Id$, or by $\delta= \Delta^2$ if $\mu\neq \Id$.
\noproof
\end{prop}

The generator $\delta$ of $Z(G)$ given in the above proposition shall be 
called the {\it standard generator} of $Z(G)$. Note also that the assumption ``$\Gamma$ is 
connected'' is not needed in (1) and (2). Let $\Gamma$ be connected. Then $\mu \neq \Id$ if and 
only if $\Gamma$ is either $A_n$, $n\ge 2$, or $D_{2n+1}$, $n\ge 2$, or $E_6$, or $I_2(2p+1)$, 
$p\ge 2$ (see \cite{BrSa}, Subsection 7.2).

Now, the following result can be found in \cite{Cha}.

\begin{prop}[Charney \cite{Cha}]
Each $a \in G$ can be uniquely written as $a=bc^{-1}$ where $b,c \in G^+$ and $b \wedge_R c =1$.
\noproof
\end{prop}

The expression $a=bc^{-1}$ of the above proposition shall be called the {\it Charney form} of $a$.

An easy observation shows that, if $s_1s_2 \dots s_l$ and $t_1 t_2 \dots t_l$ are two positive 
expressions of a same element $a \in G^+$, then the sets $\{s_1, \dots, s_l\}$ and $\{t_1, \dots, 
t_l\}$ are equal. In particular, if $a \in G_X^+$, then all the letters that appear in any 
positive expression of $a$ lie in $X$. A consequence of this fact is the following.

\begin{lemma}
Let $X$ be a subset of $S$, let $a\in G_X$, and let $a=bc^{-1}$ be the Charney form of $a$ in $G$. 
Then $b,c \in G_X^+$ and $a=bc^{-1}$ is the Charney form of $a$ in $G_X$.
\end{lemma}

\proof
Let $\vee_{X,R}$ and $\wedge_{X,R}$ denote the lattice operations of $(G_X^+,\le_R)$. The above 
observation shows that, if $a\le_R b$ and $b \in G_X^+$, then $a \in G_X^+$. This implies that $b 
\wedge_{X,R} c = b \wedge_R c$ for all $b,c \in G_X^+$. Now, let $a \in G_X$ and let $a=bc^{-1}$ 
be the Charney form of $a$ in $G_X$. We have $b,c \in G_X^+ \subset G^+$ and $b \wedge_R c = b 
\wedge_{X,R} c = 1$, thus $a=bc^{-1}$ is also the Charney form of $a$ in $G$.
\endproof

\begin{corollary}
Let $X$ be a subset of $S$. Then $G_X \cap G^+ = G_X^+$.
\noproof
\end{corollary}

\begin{corollary}
Let $X$ be a subset of $S$, $X \neq S$. Then $G_X \cap \langle \Delta \rangle = \{ 1 \}$.
\end{corollary}

\proof
Take $s \in S \setminus X$. By Proposition 2.2, we have $s \le_R \Delta$, thus $\Delta \not \in 
G_X^+ = G_X \cap G^+$.
\endproof

%%%%%%%%%%%%%%%%%%%%%%%%%%%%%%%%%%%%%%%%%%%%%%%%%%%%%%%%%%%%%%%%%%%%%%%%%%%%%%%%%%%%%%%%%%%%%
\section{Invariants}

The purpose of the present section is to calculate some invariants of the spherical type Artin 
groups.

The first invariant that we want to calculate is the cohomological dimension, denoted by $\cd 
(G)$. We assume the reader to be familiar with this notion, and we refer to \cite{Bro} for 
definitions and properties. Our result is the following.

\begin{prop}
Let $n=|S|$ be the rank of $G=G_\Gamma$. Then $\cd (G)=n$.
\end{prop}

\proof
Recall the spaces $M(W)$ and $N(W)$ defined 
in the introduction. Recall also that $\pi_1(N(W))=G$, that $W$ acts freely on $M(W)$, and that 
$N(W)=M(W)/W$. In particular, $\pi_1(M(W))$ is a subgroup of $\pi_1(N(W))=G$ (it is actually the 
kernel of the epimorphism $F: G \to W$). Finally, recall the well-known fact that, if $H_1$ 
is a subgroup of a given group $H_2$, then $\cd (H_1) \le \cd (H_2)$.

Deligne proved in \cite{Del} that $M(W)$ is aspherical, and Brieskorn proved in \cite{Bri2} that 
$H^n(M(W), \Z)$ is a free abelian group of rank $\prod_{i=1}^n m_i \neq 0$, where $m_1, m_2, \dots, m_n$ 
are the exponents of $W$, thus $n \le \cd( \pi_1( M(W)) \le \cd(G)$. On the other hand, Salvetti 
has constructed in \cite{Sal} an aspherical CW-complex of dimension $n$ whose fundamental group 
is $G$, therefore $\cd (G) \le n$.
\endproof

\bigskip
The next invariant which interests us is denoted by $\mf (G)$ and is defined to be the maximal 
order of a finite subgroup of $G/Z(G)$, where $Z(G)$ denotes the center of $G$. Its calculation 
is based on Theorems 3.2 and 3.3 given below.

Recall the permutation $\mu: S \to S$ of Proposition 2.2. This extends to an isomorphism $\mu: 
G^+ \to G^+$ which permutes the simple elements. Actually, $\mu(a)= \Delta a \Delta^{-1}$ for all 
$a \in G^+$.

\begin{thm}[Bestvina \cite{Bes}]
Assume $\Gamma$ to be connected. Let $\GG= G/ \langle \Delta^2 \rangle$, and let $H$ be a finite 
subgroup of $\GG$. Then
$H$ is a cyclic group, and,
up to conjugation, $H$ has one of the following two forms.

\smallskip\noindent
{\bf Type 1:} The order of $H$ is even, say $2p$, and
there exists a simple element $a \in T(W)$ such that $a^p=\Delta$, $\mu(a)=a$, and 
$\overline{a}$ generates $H$, where $\overline{a}$ denotes the element of $\GG$ represented by 
$a$.

\smallskip\noindent
{\bf Type 2:} The order of $H$ is odd, say $2p+1$, and
there exists a simple element $a \in T(W)$ such that $(a\, \mu(a))^{{p-1}\over 2} a = 
\Delta$ and $\overline{ a\,\mu(a)}$ generates $H$.
\noproof
\end{thm}

Now, recall the so-called {\it Coxeter number} $h$ of $W$ (see \cite{Hum}, Section 3.18). Recall 
also that this number is related to the length of $\Delta$ by the following formula
$$
{{nh}\over 2} = m_1+ \dots +m_n= \lg (\Delta)\,,
$$
where $n=|S|$ is the rank of $G$, and $m_1, \dots, m_n$ are the exponents of $W$.

\begin{thm}[Brieskorn-Saito \cite{BrSa}]
Choose any order $S=\{s_1,s_2, \dots, s_n\}$ of $S$ and write $\pi=s_1s_2 \dots s_n \in G$. Let 
$h$ be the Coxeter number of $W$.

\smallskip
(1) If $\mu = \Id$, then $h$ is even and $\pi^{h \over 2} = \Delta$.

\smallskip
(2) If $\mu \neq \Id$, then $\pi^h= \Delta^2$.
\noproof
\end{thm}

Now, we can calculate the invariant $\mf (G)$.

\begin{prop}
Assume $\Gamma$ to be connected, and let $h$ be the Coxeter number of $W$.

\smallskip
(1) If $\mu = \Id$, then $\mf (G)= h/2$.

\smallskip
(2) If $\mu \neq \Id$, then $\mf (G)=h$.
\end{prop}

\proof
Assume $\mu=\Id$. Let $\GG^0= G/Z(G) = G/ \langle \Delta \rangle$. First, observe that 
$\mf(G) \ge {h \over 2}$ by Theorem 3.3. So, it remains to prove that $\mf(G) \le {h \over 2}$, 
namely, that $|H| \le {h \over 2}$ for any finite subgroup $H$ of $\GG^0$.

Let $H$ be a finite subgroup of $\GG^0$. Consider the exact sequence
$$
1 \to \Z/ 2\Z \to \GG \stackrel{\phi}{\rightarrow} \GG^0 \to 1\,,
$$
where $\GG= G/ \langle \Delta^2 \rangle$, and set $\tilde H= \phi^{-1}(H)$. By Theorem 3.2, 
$\tilde H$ is a cyclic group and, up to conjugation, $\tilde H$ is either of Type 1 or of Type 2. 
The order of $\tilde H$ is even, say $2p$, thus $\tilde H$ is of Type 1, and there 
exists a simple element $a \in T(W)$ such that $a^p=\Delta$ and $\overline{a}$ generates $\tilde 
H$. Let $a=s_1 s_2 \dots s_r$ be an expression of $a$, and let $X=\{s_1, s_2, \dots, s_r\}$. We 
have $\Delta= a^p \in G_X$, thus, by Corollary~2.6, $X=S$ and $r=\lg (a) \ge |S|=n$. Finally,
$$
|H|={|\tilde H| \over 2} = p = {\lg (\Delta) \over \lg (a)} = \left( {nh \over 2} \right) /r \le 
{h \over 2}\,.
$$

Now, assume $\mu \neq \Id$. Let $\GG= G/Z(G)= G/\langle \Delta^2 \rangle$. First, observe that 
$\mf (G) \ge h$ by Theorem 3.3. So, it remains to prove that $\mf (G) \le h$, namely, 
that $|H| \le h$ for any finite subgroup $H$ of $\GG$.

Let $H$ be a finite subgroup of $\GG$. By Theorem 3.2, $H$ is cyclic and, up to conjugation, $H$ 
is either of Type 1 or of Type 2. Let $p$ be the order of $H$. In both cases, Type 1 and Type 2, 
there exists an element $b \in G^+$ such that $b^p=\Delta^2$ and $\overline{b}$ generates 
$H$ (take $b=a$ if $H$ is of Type 1, and $b=a\, \mu(a)$ if $H$ is of Type 2). Let $b=s_1 s_2 
\dots s_r$ be an expression of $b$, and let $X=\{s_1, s_2, \dots, s_r\}$. We have $\Delta^2=b^p 
\in G_X$, thus, by Corollary~2.6, $X=S$ and $r=\lg (b) \ge |S|=n$. It follows that
$$
|H| = p = {\lg (\Delta^2) \over \lg (b)} = {nh \over r} \le h\,.
$$
\endproof

The values of the Coxeter numbers of the irreducible Coxeter groups are well-known (see, for 
instance, \cite{Hum}, Section 3.18). Applying Proposition 3.4 to these values, one can easily 
compute the invariant $\mf (G)$ for each irreducible (spherical type) Artin group. The result is 
given in Table 1.

\def\tvia{\vrule height 0pt depth 8pt width 0pt}
\def\tvib{\vrule height 16pt depth 0pt width 0pt}
\def\tvic{\vrule height 16pt depth 8pt width 0pt}
$$\vbox{
\begin{tabular}{ccccccccccccccc}
\hline
\tvib&{\vline\kern -0.2 em \vline}&&\vline&&\vline&&\vline& $D_n,\, n\ge 4$ &\vline& $D_n,\, n\ge 
5$ &\vline&&\vline\\
\tvia$\Gamma$ &{\vline\kern -0.2 em \vline}& $A_1$ &\vline& $A_n,\ n\ge 2$ &\vline& $B_n,\, n\ge 
2$ &\vline& 
$n$ even &\vline& $n$ odd &\vline& $E_6$ &\vline\\
\hline
\tvic$\mf (G)$ &{\vline\kern -0.2 em \vline}&  1 &\vline& $n+1$ &\vline& $n$ &\vline& $n-1$ 
&\vline& $2n-2$ &\vline& 12 &\vline\\
\hline
\end{tabular}\hfill
\par
\begin{tabular}{ccccccccccccccccc}
\hline
\tvib&{\vline\kern -0.2 em \vline}&&\vline&&\vline&&\vline&&\vline&&\vline& $I_2(p),\, p\ge 6$ 
&\vline& $I_2(p),\, p\ge 5$ &\vline\\
\tvia$\Gamma$ &{\vline\kern -0.2 em \vline}& $E_7$ &\vline& $E_8$ &\vline& $F_4$ &\vline& 
$H_3$ &\vline& $H_4$ &\vline& $p$ even &\vline& $p$ odd &\vline\\
\hline
\tvic$\mf (G)$ &{\vline\kern -0.2 em \vline}&  9 &\vline& 15 &\vline& 6 &\vline& 5 &\vline& 15 
&\vline& $p/2$ &\vline& $p$ &\vline\\
\hline
\end{tabular}\hfill}$$
\bigskip
\centerline{{\bf Table 1:} The invariant $\mf(G)$.}

\bigskip\noindent
{\bf Remark.} Combining \cite{Bes}, Theorem 4.5, with \cite{BrMi}, Section 3, one can actually 
compute all the possible orders for a finite subgroup of $G/Z(G)$. The maximal order suffices for 
our purpose, thus we do not include this more complicate calculation in this paper.

\bigskip
The next invariant that we want to compute is the rank of the abelianization of $G$ that we denote 
by $\rkAb (G)$. This invariant can be easily computed using the standard presentation of $G$, and 
the result is as follows.

\begin{prop}
Let $\Gamma_0$ be the (non-labelled) graph defined by the following data.

\smallskip
$\bullet$ $S$ is the set of vertices of $\Gamma_0$;

\smallskip
$\bullet$ two vertices $s,t$ are joined by an edge if $m_{s\,t}$ is odd.

\smallskip\noindent
Then the abelianization of $G$ is a free abelian group of rank $\rkAb (G)$, the number of connected 
components of $\Gamma_0$.
\noproof
\end{prop}

The last invariant which interests us is the rank of the center of $G$ that we denote by $\rkZ 
(G)$. The following proposition is a straightforward consequence of Proposition 2.2.

\begin{prop}
The center of $G$ is a free abelian group of rank $\rkZ (G)$, the number of components of $\Gamma$.
\noproof
\end{prop}

%%%%%%%%%%%%%%%%%%%%%%%%%%%%%%%%%%%%%%%%%%%%%%%%%%%%%%%%%%%%%%%%%%%%%%%%%%%%%%%%%%%%%%%%%%%%%%%%%
\section{Irreducibility}

Throughout this section, we assume that $G$ is irreducible (namely, that $\Gamma$ is connected). 
Let $H_1, H_2$ be two subgroups of $G$. Recall that $[H_1, H_2]$ denotes the subgroup of $G$ 
generated by $\{a_1^{-1} a_2^{-1} a_1a_2; a_1 \in H_1\text{ and } a_2 \in H_2\}$. Our goal in this 
section is to show that $G$ cannot be expressed as $G= H_1 \cdot H_2$ with $[H_1, H_2] = \{1\}$, 
unless either $H_1 \subset Z(G)$ or $H_2 \subset Z(G)$. This shall implies that $G$ cannot be a 
non-trivial direct product.

Recall that $\delta$ denotes the standard generator of $Z(G)$. For $X \subset S$, we denote by 
$\delta_X$ the standard generator of $G_X$, and, for $a \in G$, we denote by $Z_G(a)$ the 
centralizer of $a$ in $G$.

\begin{lemma}
Let $t \in S$ such that $\Gamma_{S \setminus \{t\}}$ is connected and $\mu(t) \neq t$ if $\mu 
\neq \Id$. Then $Z_G( \delta_{S \setminus \{t\}})$ is generated by $G_{S \setminus \{t\}} \cup \{ 
\delta \}$ and is isomorphic to $G_{S \setminus \{t\}} \times \{ \delta \}$.
\end{lemma}

\proof
Assume first that $\mu=\Id$ (in particular, $\delta=\Delta$). By \cite{Par1}, Theorem 5.2, 
$Z_G(\delta_{S\setminus \{t\}})$ is generated by $G_{S \setminus \{t\}} \cup \{ \Delta^2, \Delta 
\Delta_{S \setminus \{t\}}^{-1}\}$, thus $Z_G(\delta_{S \setminus \{t\}})$ is generated by $G_{S 
\setminus \{t\}} \cup \{ \delta \}$.

Now, assume $\mu \neq \Id$ (in particular, $\delta = \Delta^2$ and $\mu(t) \neq t$). By 
\cite{Par1}, Theorem 5.2, $Z_G(\delta_{S \setminus \{t\}})$ is generated by $G_{S \setminus 
\{t\}} \cup \{ \Delta^2, \Delta \Delta_{S \setminus \{\mu(t)\}}^{-1} \Delta \Delta_{S \setminus 
\{t\}}^{-1} \}$. Observe that 
\break
$\Delta \Delta_{S \setminus \{\mu(t)\}}^{-1} \Delta \Delta_{S 
\setminus \{t\}}^{-1} = \Delta^2 \Delta_{S \setminus \{t\}}^{-2}$, thus $Z_G(\delta_{S \setminus 
\{t\}})$ is generated by $G_{S \setminus \{t\}} \cup \{\delta\}$.

By the above, we have an epimorphism $G_{S \setminus \{t\}} \times \langle \delta \rangle \to 
Z_G( \delta_{S \setminus \{t\}})$, and, by Corollary~2.6, the kernel of this epimorphism is 
$\{1\}$.
\endproof

\bigskip\noindent
{\bf Remark.} It is an easy exercise to show (under the assumption that $\Gamma$ is connect) that 
there always exists $t \in S$ such that $\Gamma_{S \setminus \{t\}}$ is connected and $\mu(t) 
\neq t$ if $\mu \neq \Id$.

\begin{prop}
Let $H_1, H_2$ be two subgroups of $G$ such that $G=H_1 \cdot H_2$ and $[H_1,H_2]=\{1\}$. Then 
either $H_1 \subset Z(G)$ or $H_2 \subset Z(G)$. If, moreover, $H_1 \cap H_2 = \{1\}$, then 
either $H_1=\{1\}$ and $H_2=G$, or $H_1=G$ and $H_2= \{1\}$.
\end{prop}

\proof
We argue by induction on $n= |S|$. If $n=1$, then $\Gamma= A_1$ and $G=\Z$, and the conclusion of 
the proposition is well-known.

Assume $n \ge 2$. For $i=1,2$, let $\tilde H_i$ denote the subgroup of $G$ generated by $H_i \cup 
\{ \delta\}$. We have $G= \tilde H_1 \cdot \tilde H_2$, $[\tilde H_1, \tilde H_2] = \{1\}$, $H_1 
\subset \tilde H_1$, and $H_2 \subset \tilde H_2$. Observe also that $\tilde H_1 \cap \tilde H_2$ 
must be included in the center of $G$, and that $\delta \in \tilde H_1 \cap \tilde H_2$, thus 
$\tilde H_1 \cap \tilde H_2 = \langle \delta \rangle$. Take $t \in S$ such that $\Gamma_{S 
\setminus \{t\}}$ is connected and $\mu(t) \neq t$ if $\mu \neq \Id$, write $X=S \setminus 
\{t\}$, and choose $d_1 \in \tilde H_1$ and $d_2\in \tilde H_2$ such that $\delta_X=d_1d_2$.

Let $a \in G_X$. Choose $a_1 \in \tilde H_1$ and $a_2 \in \tilde H_2$ such that $a=a_1a_2$. We 
have 
$$
1= a^{-1} \delta_X^{-1} a \delta_X = a_1^{-1} d_1^{-1} a_1d_1 a_2^{-1} d_2^{-1} a_2d_2\,,
$$
thus
$$
a_1^{-1} d_1^{-1} a_1d_1 = d_2^{-1} a_2^{-1} d_2a_2 \in \tilde H_1 \cap \tilde H_2 = \langle 
\delta \rangle \,.
$$
Let $k \in \Z$ such that $a_1^{-1} d_1^{-1} a_1d_1 = \delta^k$. Consider the homomorphism $\deg: 
G \to \Z$ which sends $s$ to $1$ for all $s \in S$. Then
$$
0= \deg( a_1^{-1} d_1^{-1} a_1d_1) = \deg( \delta^k)= k\,\lg (\delta)\,,
$$
thus $k=0$, hence $a_1$ and $d_1$ commute. Now, $a_1$ and $d_2$ also commute (since $a_1 \in 
\tilde H_1$ and $d_2 \in \tilde H_2$), thus $a_1$ commutes with $\delta_X = d_1d_2$. By Lemma 
4.1, $a_1$ can be written as $a_1= b_1 \delta^{p_1}$, where $b_1 \in G_X$ and $p_1 \in \Z$. Note 
also that $b_1 = a_1 \delta^{-p_1} \in \tilde H_1$, since $\delta \in \tilde H_1$, thus 
$b_1 \in G_X \cap \tilde H_1$. 
Similarly, $a_2$ 
can be written as $a_2=b_2 \delta^{p_2}$ where $b_2 \in G_X \cap \tilde H_2$ and $p_2 \in \Z$. We 
have $\delta^{p_1+p_2} = a b_1^{-1} b_2^{-1} \in G_X \cap \langle \delta \rangle = \{1\}$ (by 
Corollary 2.6), thus $p_1+p_2=0$ and $a=b_1b_2$.

So, we have 
$$
G_X= (G_X \cap \tilde H_1) \cdot (G_X \cap \tilde H_2)\,.
$$
Moreover, by Corollary 2.6,
$$
(G_X \cap \tilde H_1) \cap (G_X \cap \tilde H_2) = G_X \cap \langle \delta \rangle = \{1\}\,.
$$
By the inductive hypothesis, it follows that, up to permutation of 1 and 2, we 
have $G_X \cap \tilde H_1 = G_X$ (namely, $G_X \subset \tilde H_1$), and $G_X \cap \tilde H_2 = 
\{1\}$.

We turn now to show that $\tilde H_2 \subset \langle \delta \rangle = Z(G)$. Since $H_2 \subset 
\tilde H_2$, this shows that $H_2 \subset Z(G)$.

Let $a \in \tilde H_2$. Since $\delta_X \in G_X \subset \tilde H_1$, $a$ and $\delta_X$ commute. 
By Lemma 4.1, $a$ can be written as $a=b \delta^p$, where $b \in G_X$ and $p \in \Z$. Since 
$\delta \in \tilde H_2$, we also have $b = a \delta^{-p} \in \tilde H_2$, thus $b \in G_X \cap 
\tilde H_2 = \{1\}$, therefore $a = \delta^p \in \langle \delta \rangle$.

Now, assume that $H_1 \cap H_2=\{1\}$. By the above, we may suppose that $H_2 \subset Z(G)= 
\langle \delta \rangle$. In particular, there exists $k \in \Z$ such that $H_2 = \langle 
\delta^k \rangle$.
Choose any order $S=\{s_1, \dots, s_n\}$ of $S$, and write $\pi=s_1s_2 \dots s_n \in G$.
Let $b \in H_1$ and $p \in \Z$ such that $\pi= b \delta^{pk}$.
Observe that $b \neq 1$ since $\pi$ is not central in $G$. Let $h$ denote the Coxeter number of 
$W$. By Theorem~3.3, $\pi^h= b^h \delta^{phk} \in Z(G)$, thus $b^h \in Z(G)$. Moreover, $b^h \neq 
1$ since $G$ is torsion free and $b \neq 1$. 
This implies that $Z(H_1)\neq\{1\}$. Now, observe that $Z(H_1) \subset Z(G)= \langle 
\delta \rangle$, thus there exists $l >0$ such that $Z(H_1)=\langle \delta^l \rangle$. Finally, 
$\delta^{lk} \in H_1 \cap H_2 = \{1\}$, thus $kl=0$, therefore $k=0$ (since $l\neq 0$) and 
$H_2=\{1\}$. Then we also have $H_1=G$.
\endproof

\begin{prop}
Assume $n=|S| \ge 2$. Let $H$ be a subgroup of $G$ such that $G=H \cdot \langle \delta \rangle$. 
Then $\cd (H)= \cd (G)$, $\mf(H) = \mf(G)$, and $\rkAb (H)= \rkAb(G)$.
\end{prop}

\proof
For all $s \in S$, take $b_s \in H$ and $p_s \in \Z$ such that $s=b_s \delta^{p_s}$. We can and 
do suppose that $p_s=p_t$ if $s$ and $t$ are conjugate in $G$. Then the mapping $S \to H$, $s 
\mapsto b_s=s\delta^{-p_s}$ determines a homomorphism $\varphi: G \to H$.

We show that $\varphi: G \to H$ is injective. Observe that the mapping $S \to \Z$, $s \mapsto 
p_s$ determines a homomorphism $\eta: G \to \Z$, and that $\varphi(a)= a \delta^{-\eta (a)}$ for 
all $a \in G$. In particular, if $a \in \Ker \varphi$, then $a=\delta^{\eta(a)} \in Z(G)$. Choose 
any order $S=\{s_1, \dots, s_n\}$ of $S$, and write $\pi=s_1s_2 \dots s_n \in G$. Note that 
$\varphi(\pi) \neq 1$, since $\pi$ is not central in $G$, and that, by Theorem 3.3, there exists 
$k>0$ such that $\pi^k = \delta$. Let $a \in \Ker \varphi$. Then $a = \delta^{\eta (a)} = \pi^{k 
\eta(a)}$, thus $1=\varphi(a)= \varphi(\pi)^{k \eta(a)}$. We have $\varphi(\pi) \neq 1$ and $G$ 
is torsion free, hence $\eta(a) =0$ (since $k>0$) and $a=1$.

Now, recall that $\cd (H_1) \le \cd (H_2)$ if $H_1$ is a subgroup of a given group $H_2$. So,
$$
\cd (G)= \cd (\varphi(G)) \le \cd(H) \le \cd (G)\,.
$$

The equality $G=H \cdot \langle \delta \rangle = H \cdot Z(G)$ implies that $Z(H)= Z(G) \cap H$ 
and $G/Z(G)= H/Z(H)$. In particular, we have $\mf(H)= \mf(G)$.

Let $\HH$ be a group, let $g$ be a central element in $\HH$, and let $p>0$. Let $\GG= (\HH \times \Z)/ 
\langle (g,p) \rangle$. Then one can easily verify (using the Reidemeister-Schreier method, for 
example) that we have exact sequences $1 \to \HH \to \GG \to \Z/p\Z \to 1$ and $1 \to \Ab (\HH) 
\to \Ab (\GG) \to \Z/p\Z \to 1$, where $\Ab (\GG)$ (resp. $\Ab(\HH)$) denotes the abelianization 
of $\GG$ (resp. $\HH$).

Now, recall the equality $G=H \cdot \langle \delta \rangle$. By Proposition 4.2, we have $H \cap 
\langle \delta \rangle \neq \{1\}$. So, there exists $p>0$ such that $H \cap \langle \delta 
\rangle = \langle \delta^p \rangle$. Write $d=\delta^p \in H$. Then $d$ is central in $H$ and $G 
\simeq (H \times \Z)/\langle (d,p) \rangle$. By the above observation, it follows that we have an 
exact sequence $1 \to \Ab (H) \to \Ab (G) \to \Z/p\Z \to 1$, thus $\Ab (H)$ is a free abelian 
group of rank $\rkAb (G)$.
\endproof
 
%%%%%%%%%%%%%%%%%%%%%%%%%%%%%%%%%%%%%%%%%%%%%%%%%%%%%%%%%%%%%%%%%%%%%%%%%%%%%%%%%%%%%%%%%%%%%%
\section{Proof of the main theorem}

\begin{prop}
Let $\Gamma$ and $\Omega$ be two connected spherical type Coxeter graphs, and let $G$ and $H$ be 
the Artin groups associated to $\Gamma$ and $\Omega$, respectively. If $\cd (G)= \cd (H)$, $\mf 
(G)= \mf (H)$, and $\rkAb(G)= \rkAb(H)$, then $\Gamma= \Omega$.
\end{prop}

\proof
Let $n$ and $m$ be the numbers of vertices of $\Gamma$ and $\Omega$, respectively. By Proposition~3.1, 
we have $n= \cd(G)= \cd(H) = m$.

Suppose $n=m=1$. Then $\Gamma= \Omega= A_1$.

Suppose $n=m \ge 3$. Then one can easily verify in Table 1 that the equality $\mf (G)= \mf (H)$ 
implies $\Gamma= \Omega$.

Suppose $n=m=2$. Let $p,q \ge 3$, such that $\Gamma= I_2(p)$ and $\Omega= I_2(q)$. By Proposition~3.5, 
either $\rkAb(G)=\rkAb(H)=2$ and $p,q$ are both even, or $\rkAb(G)= \rkAb(H)=1$ and $p,q$ 
are both odd. If $p,q$ are both even, then, by Table 1, ${p \over 2} = \mf(G) = \mf(H)= {q 
\over 2}$, thus $p=q$ and $\Gamma=\Omega= I_2(p)$. If $p,q$ are both odd, then, by Table 1, $p= 
\mf (G)= \mf (H)=q$, thus $\Gamma=\Omega= I_2(p)$.
\endproof

\begin{corollary}
Let $\Gamma$ and $\Omega$ be two connected spherical type Coxeter graphs, and let $G$ and $H$ be 
the Artin groups associated to $\Gamma$ and $\Omega$, respectively. If $G$ is isomorphic to $H$, 
then $\Gamma=\Omega$.
\noproof
\end{corollary}

\noindent
{\bf Proof of Theorem 1.1.} 
Let $\Gamma$ and $\Omega$ be two spherical type Coxeter graphs, and let $G$ and $H$ be the Artin 
groups associated to $\Gamma$ and $\Omega$, respectively. We assume that $G$ is isomorphic to $H$ 
and turn to prove that $\Gamma=\Omega$.

Let $\Gamma_1, \dots, \Gamma_p$ be the connected components of $\Gamma$, and let $\Omega_1, 
\dots, \Omega_q$ be the connected components of $\Omega$. For $i=1, \dots, p$, we denote by $G_i$ 
the Artin group associated to $\Gamma_i$, and, for $j=1, \dots, q$, we denote by $H_j$ the Artin 
group associated to $\Omega_j$. We have $G= G_1 \times G_2 \times \dots \times G_p$ and $H=H_1 
\times H_2 \times \dots \times H_q$. We may and do assume that there exists $x \in \{0,1, \dots, 
p\}$ such that $\Gamma_i \neq A_1$ for $i=1, \dots, x$, and $\Gamma_i=A_1$ for $i=x+1, \dots, p$. 
So, $G_1, \dots, G_x$ are non abelian irreducible Artin groups of rank $\ge 2$, and $G_{x+1}, 
\dots, G_p$ are all isomorphic to $\Z$. Similarly, we may and do assume that there exists $y \in 
\{0,1, \dots, q\}$ such that $\Omega_j \neq A_1$ for $j=1, \dots, y$, and $\Omega_j=A_1$ for 
$j=y+1, \dots, q$. We can also assume that $x \ge y$.

A first observation is, by Proposition 3.6, that
$$
p= \rkZ(G)= \rkZ(H)= q\,.
$$

Now, fix an isomorphism $\varphi: G \to H$. For $1 \le i\le p$, let $\iota_i: G_i \to G$ be the 
natural embedding, for $1 \le j\le p$, let $\kappa_j: H \to H_j$ be the projection on the $j$-th 
component, and, for $1 \le i,j \le p$, let $\varphi_{i\,j} = \kappa_j \circ \varphi \circ 
\iota_i: G_i \to H_j$.

Let $j\in \{1, \dots, y\}$. Observe that $H_j= \prod_{i=1}^p \varphi_{i\,j}(G_i)$, and that 
$[\varphi_{i\,j}(G_i), \varphi_{k\,j}(G_k)]=1$ for all $i,k \in \{1, \dots, p\}$, $i \neq k$. Let 
$\delta_j^H$ denote the standard generator of $Z(H_j)$, and, for $i\in\{1, \dots, p\}$, let $\tilde 
H_{i\,j}$ be the subgroup of $H_j$ generated by $\varphi_{i\,j}(G_i) \cup \{\delta_j^H\}$. By 
Proposition~4.2, there exists $\chi (j) \in \{1, \dots, p\}$ such that $H_j= \tilde 
H_{\chi(j)\,j}$, and $\tilde H_{i\,j} = Z(H_j)= \langle \delta_j^H \rangle$ for $i \neq \chi(j)$. 
Since $H_j$ is non abelian, $\chi(j)$ is unique and $\chi(j)\in \{1, \dots, x\}$.

We turn now to show that the map $\chi: \{1, \dots, y\} \to \{1, \dots, x\}$ is surjective. Since 
$x \ge y$, it follows that $x=y$ and $\chi$ is a permutation.

Let $i \in \{1, \dots, x\}$ such that $\chi(j) \neq i$ for all $j \in \{1, \dots, y\}$. Then 
$\varphi_{i\,j} (G_i) \subset Z(H_j)$ for all $j=1, \dots, p$, thus $\varphi(G_i) \subset Z(H)$. 
This contradicts the fact that $\varphi$ is injective and $G_i$ is non abelian.

So, up to renumbering the $\Gamma_i$'s, we can suppose that $\chi(i)=i$ for all $i\in \{1, \dots, 
x\}$.

We prove now that $\varphi_{i\,i}: G_i \to H_i$ is injective for all $i \in \{1, \dots, x\}$. Let 
$a \in \Ker \varphi_{i\,i}$. Since $\varphi_{i\,j} (a) \in Z(H_j)$ for all $j \neq i$, we have 
$\varphi(a) \in Z(H)$. Since $\varphi$ is injective, it follows that $a \in Z(G_i)$. Let $\{s_1, 
\dots, s_r\}$ be the set of vertices of $\Gamma_i$, and let $\pi= s_1s_2 \dots s_r \in G_i$. 
Observe that $\varphi_{i\,i}(\pi) \neq 1$ since $\pi$ is not central in $G_i$. Let $\delta_i^G$ 
be the standard generator of $Z(G_i)$. By Theorem 3.3, there exists $k >0$ such that $\pi^k= 
\delta_i^G$. On the other hand, since $a \in Z(G_i)$, there exists $l \in \Z$ such that $a= 
(\delta_i^G)^l= \pi^{kl}$. Now, $1= \varphi_{i\,i}(a)= \varphi_{i\,i}(\pi)^{kl}$, $H_i$ is 
torsion free, and $\varphi_{i\,i}(\pi) \neq 1$, thus $kl=0$ and $a=\pi^{kl}=1$.

Let $i \in \{1, \dots, x\}$. Recall that $\varphi_{i\,i}: G_i \to H_i$ is injective, and $H_i= 
\varphi_{i\,i} (G_i) \cdot \langle \delta_i^H \rangle$, where $\delta_i^H$ denotes the standard 
generator of $H_i$. By Proposition 4.3, it follows that
$$
\cd(G_i)= \cd(H_i)\,,\quad \mf(G_i) =\mf(H_i)\,,\quad \rkAb(G_i)= \rkAb(H_i)\,,
$$
thus, by Proposition 5.1, $\Gamma_i=\Omega_i$. Let $i \in \{x+1, \dots, p\}$. Then $\Gamma_i= 
\Omega_i = A_1$. So, $\Gamma= \Omega$.
\endproof

\bigskip\noindent
{\bf Remark.} In the proof above, the homomorphism $\varphi_{i\,i}$ is injective but is not 
necessarily surjective as we show in the following example.

Let $G_1= \langle s_1,s_2 | s_1s_2s_1= s_2s_1s_2 \rangle$ be the Artin group associated to $A_2$, 
let $G_2=\Z= \langle t \rangle$, and let $G=G_1 \times G_2$. We denote by $\delta= (s_1s_2)^3$ 
the standard generator of $Z(G_1)$. Let $\varphi: G \to G$ be the homomorphism defined by
$$
\varphi(s_1)= s_1 \delta t\,, \quad \varphi(s_2)= s_2 \delta t \,, \quad \varphi(t)= \delta t \,.
$$
Then $\varphi$ is an isomorphism but $\varphi_{1\,1}$ is not surjective. The inverse $\varphi^{-
1}: G \to G$ is determined by
$$
\varphi^{-1} (s_1)= s_1 t^{-1}\,, \quad \varphi^{-1} (s_2)= s_2 t^{-1}\,, \quad \varphi^{-1} (t)= 
\delta^{-1} t^7\,.
$$

%%%%%%%%%%%%%%%%%%%%%%%%%%%%%%%%%%%%%%%%%%%%%%%%%%%%%%%%%%%%%%%%%%%%%%%%%%%%%%%%%%%%%%%%

%%%%%%%%%%%%%%%%%%%%%%%%%%%%%%%%%%%%%%%%%%%%%%%%%%%%%%%%%%%%%%%%%%%%%%%%%%%%%%%%%%%%%%%%%

\bigskip\bigskip\noindent
\halign{#\hfill\cr
Luis Paris\cr
Institut de Math\'ematiques de Bourgogne\cr
Universit\'e de Bourgogne\cr
UMR 5584 du CNRS, BP 47870\cr
21078 Dijon cedex\cr
FRANCE\cr
\noalign{\smallskip}
\texttt{lparis@u-bourgogne.fr}\cr}

\end{document}